\documentclass[final,leqno]{siamltex}
\usepackage{amsmath,bm}
\usepackage{amsfonts}
\usepackage{color}
\usepackage[latin1]{inputenc}
\usepackage{algpseudocode}
\usepackage{algorithm}
\usepackage{multirow}
\usepackage{caption}
\usepackage{mathrsfs}
\usepackage{lscape}
\usepackage{mwe}
\usepackage{subfig}
\usepackage{url}
\usepackage[update,prepend]{epstopdf} 
\usepackage{graphicx}

\newcommand{\be}{\begin{eqnarray}}
\newcommand{\e}{\end{eqnarray}}
\newcommand{\bes}{\begin{eqnarray*}}
\newcommand{\es}{\end{eqnarray*}}
\newcommand{\beq}{\begin{equation}}
\newcommand{\eeq}{\end{equation}}

\newtheorem{remark}{Remark}
\author{H. C.  Elman and A. Onwunta}

\usepackage[pagewise]{lineno}

\begin{document}
\title{Reduced-order modeling for  nonlinear Bayesian statistical inverse 
problems\thanks{This work was supported by the U.S. Department of Energy Office of Advanced Scientific Computing Research,
Applied Mathematics program under award DE-SC0009301 and by the U.S. National Science Foundation under grant DMS1819115.}}
\author{Howard C. Elman\footnotemark[1] and Akwum Onwunta\footnotemark[2]}
\renewcommand{\thefootnote}{\fnsymbol{footnote}}
\footnotetext[1]{Department of Computer Science and Institute for Advanced Computer Studies,
University of Maryland, College Park,
MD 20742, USA,
(\url{elman@cs.umd.edu}).}
\footnotetext[2]{Department of Computer Science,
University of Maryland, College Park,
MD 20742, USA,
(\url{onwunta@cs.umd.edu}).}
 \renewcommand{\thefootnote}{\arabic{footnote}}
\maketitle
\begin{abstract}
Bayesian statistical inverse problems 
are often  solved with Markov chain Monte Carlo (MCMC)-type schemes.
When the problems are governed by large-scale discrete nonlinear partial differential
equations (PDEs), they  are computationally challenging because one would then need to solve the forward problem at every
sample point. In this paper, the use of the discrete empirical interpolation  method (DEIM)  is considered
for the forward solves within an  MCMC routine. A preconditioning
strategy for the DEIM model   is also proposed to 
accelerate the forward solves. The preconditioned DEIM model 
is applied to a finite difference discretization of  a nonlinear PDE
in the MCMC model. Numerical experiments
show that this approach yields accurate forward results. Moreover, the computational cost 
of solving the associated statistical inverse problem is reduced by more than $70\%$.
\end{abstract}


\pagestyle{myheadings}
\thispagestyle{plain}
\markboth{Efficient solver for nonlinear statistical inverse problems}{Efficient reduced-order modeling for nonlinear statistical inverse problems}
\section{Introduction}
\label{sec1}
In the last few decades, computational science and engineering has seen a dramatic increase in the need for 
simulations 
of large-scale  inverse problems \cite{BSHL2014, Bui14, CS07, Flath11, KS05}.  
When these problems are governed by partial differential equations (PDEs),
the ultimate goal is to recover quantities from limited and noisy observations. 
One of the major challenges posed by such problems is that they are  computationally
expensive and, as result, entail vast storage requirements. Moreover, they are often ill-posed -- a feature
which can induce uncertainty in 
the recovered quantities \cite{Flath11}. 
Generally speaking, solving an inverse problem could be accomplished using either a deterministic technique
 or  a Bayesian statistical method. The 
  deterministic approach essentially leads to an optimization problem in which the objective
  function is minimized to obtain a point estimate for the sought parameter. 
However, this method does not take into consideration possible uncertainties in the solution of
the inverse problem. Instead of computing a single solution to the inverse problem, 
the Bayesian inference framework, on the other hand, seeks  a probability distribution
of a set of plausible solutions. More precisely, the Bayesian statistical approach is a systematic
way of modeling the  observed data and parameters of interest as random variables to enable
the quantification of the uncertainties in the reconstructed quantities \cite{KS05}. The outcome 
is the so-called posterior distribution as the solution of the inverse problem, from which one 
can determine the mean, covariance, and other statistics of the unknown parameters.

The appealing features of the Bayesian statistical method notwithstanding,
the resulting posterior
distribution quite often does not admit analytic expression.
Consequently, efficient numerical methods 
are required to approximate the solution to the inverse problem. In computational practice,  
 the Markov chain Monte
Carlo (MCMC) sampling techniques are probably the commonest methods used to explore the posterior distribution
 \cite{BSHL2014, Bui14, CS07, KS05, Martin12, RCS14}. For  large-scale inverse problems involving 
 discretizations of nonlinear PDEs, however, the MCMC schemes can be
 computationally intractable because they require the solution of 
 high-dimensional discrete nonlinear PDEs at every
sample point.
This is  the main focus of this work; here, we employ 
state-of-the-art nonlinear model order reduction techniques within  MCMC algorithms 
to  reduce the computational complexity of the statistical inverse problem.
More specifically, following \cite{BMNP04, CS10, EF17, Forstall2015} and the
references therein, we apply the 
POD-Galerkin method and a {\it preconditioned} discrete empirical interpolation  method (DEIM) to 
accelerate the forward solves in our MCMC routine.

The outline of the paper is  as follows. A general formulation of the discrete 
nonlinear PDE model, as well as the associated statistical inverse problem we would like to solve, is 
presented in Section \ref{prob_form}. The
MCMC scheme for sampling the posterior density of the unknown parameters in the inverse 
problem is presented in Section \ref{mcmc}.
Section \ref{rom_deim} discusses the reduced-order models which will be employed to facilitate the forward solves
within the MCMC algorithm. Finally,  Section \ref{numexp}  reports our preconditioning strategy for the forward 
solves, as well as the numerical
experiments.

\section{Problem formulation}
\label{prob_form}
In this section, we present our nonlinear statistical inverse problem
 of interest. To begin with, we present a general
 formulation of the governing discrete nonlinear PDE  model in Section 
 \ref{nonl}. Next, we
 proceed to discuss in Section \ref{sinp}
 the associated Bayesian statistical inverse problem.  
\subsection{Parameterized  nonlinear full model}
\label{nonl}
 We consider
here a system of parameterized  nonlinear equations that
 results from spatial discretization of a nonlinear PDE using, for instance, a finite difference or
finite element method. More precisely, we are interested in the following forward model

\be
\label{disc1}
G(u;\xi) := Au + F(u;\xi) + B =0, 
\e
where the vectors $u(\xi)\in \mathbb{R}^N$  and  
 $\xi=[\xi_1, \dots, \xi_d]^T\in \mathcal{D},$ denote  the  state and the  parameters
in parameter domain $\mathcal{D}\subset\mathbb{R}^d,$ respectively. Moreover, the  matrix
$A\in\mathbb{R}^{N\times N}$  
represents the linear operators in the  governing PDE, and  $B$  a source term. 
The function $F:\mathbb{R}^N\times \mathcal{D} \rightarrow \mathbb{R}^N$ corresponds 
to all the nonlinear terms in the  PDE.
The quantity of interest $y(\xi)$ associated with the forward model (\ref{disc1})  is 
expressed as
\be
\label{output}
y(\xi) = Cu(\xi), 
\e
where 
$C\in\mathbb{R}^{m\times N}$  is also a constant  matrix.

In what follows, our ultimate goal is  the case where the  state vector $u$  is high-dimensional,
and the numerical solution of the nonlinear system (\ref{disc1})  is
needed at many parameter values.

\subsection{Bayesian statistical inverse problem}
\label{sinp}
In what follows, we  present the Bayesian formulation and solution method for our statistical inverse model problem.
We refer to \cite{CS07, KS05} for  excellent introductions to statistical inverse problems.
Note that the goal of the  inverse
problem associated with  the forward problem (\ref{disc1}) is essentially 
to estimate  a parameter vector $\xi \in\mathcal{D}$ given some observed data
$y_{obs}\in\mathbb{R}^m.$ 
To achieve this goal,  the Bayesian  framework assumes that both the parameters 
of interest $\xi$ and the observed data  $y_{obs}$ are random variables.
Moreover, the solution to
statistical inverse problem is a posterior probability density,
$\pi(\xi|y_{obs}):\mathbb{R}^d\rightarrow \mathbb{R},$ 
which encodes the uncertainty  from the set of observed data and the sought parameter vector.
Mathematically, $\pi(\xi|y_{obs})$ is described by two density functions: the prior
$\pi(\xi):\mathbb{R}^d\rightarrow \mathbb{R}$   
and the likelihood function $\pi(y_{obs}|\xi)$.  
More precisely,  the  prior probability density $\pi(\xi)$ 
describes the information  which we would like to 
enforce on the parameters $\xi$  before considering the observed data, while
the likelihood function $\pi(y_{obs}|\xi)$ is  a conditional
probability that models the relationship
between  the observed data $y_{obs}$ and 
 set of unknown parameters $\xi$.

From
Bayes' theorem, the posterior probability density  is given by
\begin{equation}
\label{postnew}
\pi(\xi|y_{obs})
     =   \frac{\pi(\xi) \pi(y_{obs}|\xi)}{  \pi(y_{obs}) }
     \propto  \pi(\xi)  \pi(y_{obs}|\xi),
\end{equation}
where the
quantity $\pi(y_{obs})$ is a normalization constant.
We follow \cite{CS07, Flath11} to derive a computable expression for the posterior 
probability density $\pi(\xi|y_{obs})$. To this end, we assume that the observed data  are of the form
\be
\label{obs}
 y_{obs}=y(\xi) + \varepsilon,
\e
where $\xi$ is the vector of unknown parameters that we wish to estimate,
$y$ is the so-called nonlinear {\it parameter-to-observable} map, 
and $\varepsilon\in\mathbb{R}^m$ is an $m$-dimensional Gaussian error with zero mean and covariance
$\Sigma_{\varepsilon}\in\mathbb{R}^{m\times m}$, with $\Sigma_{\varepsilon} = \sigma^2 I.$
That is, $\varepsilon \sim \mathcal{N}(0,\sigma^2 I)$. Under the assumption
of statistical independence between the noise $\varepsilon$ and the unknown parameters $\xi,$  it follows that 
 the likelihood function $\pi(y_{obs}|\xi) =\pi(y_{obs} -y(\xi)) $ so that (see, e.g., \cite[p. 42]{CS07})
\be
\label{like}
\pi(y_{obs}|\xi)
     & \propto &  \exp\left[  -\frac{1}{2\sigma^2} (y_{obs}-y(\xi))^T (y_{obs}-y(\xi))  \right].
\e

Similar to the likelihood function, the prior for $\xi$ can also be modeled as Gaussian \cite{Flath11}. In this work, however, we assume,
for simplicity, that $\xi$
is uniformly distributed over $\mathcal{D}$. Consequently, 
the posterior probability density in (\ref{postnew})  is expressed as
\begin{align}
\label{postf}
 \pi(\xi|y_{obs}) &\propto
  \begin{cases}
    \exp\left[  -\frac{1}{2\sigma^2} (y_{obs}-y(\xi))^T (y_{obs}-y(\xi))  \right]        & \text{if } \xi \in \mathcal{D}, \\
   0        & \text{otherwise.}
  \end{cases}
\end{align}
Observe from the expression (\ref{postf}) that due to the nonlinearity of $y,$ it is
quite difficult, if not impossible, to sample directly from $\pi(\xi|y_{obs})$ in (\ref{postf}).
It is a common computational practice to use Markov chain Monte Carlo algorithms to 
tackle this problem \cite{Bui14, Martin12}.

\section{Markov chain Monte Carlo techniques}
\label{mcmc}

In this ssection, we discuss the evaluation of the  posterior probability density
 as defined in (\ref{postf}). First,  note that the evaluation of 
 the  posterior probability density could be accomplished via quadrature techniques
 or  Monte Carlo integration techniques. However, these methods can be quite computationally intensive 
 when one is dealing with high-dimensional problems \cite{RCS14}. An alternative method,
 which we employ in this work, is the Markov chain Monte Carlo (MCMC) technique
 \cite{BSHL2014, Bui14, CS07, RCS14}.
 Unlike the quadrature 
 and the  Monte Carlo integration techniques, the MCMC method uses the attributes of
 the density to specify parameter values that adequately explore the geometry 
 of the distribution. This is achieved by constructing Markov chains whose stationary
 distribution is the posterior density.

There are different ways to implement MCMC; the prominent among them include
Metropolis-Hastings (MH) \cite{RCS14}, Gibbs \cite{BSHL2014} and Hamiltonian \cite{Bui14}
 algorithms.
The MH  algorithm  can be conceptualized as a random walk through the
parameter space that  starts at an arbitrary point and where the
next step  depends solely on the
current position (this is  the Markov property of the MCMC). 
Moreover,  each random sample is generated from a {\it proposal distribution} (a term that will be 
made precise in the next section) which is subject to very mild restrictions.
The Gibbs algorithm can be thought of as a special case of the MH 
algorithm in which the proposal distribution is determined by both the component parameter selected
and its current location in the parameter space. 
Similar to the Gibbs algorithm,
the Hamiltonian  algorithm  employs a dynamical  proposal distribution.
However, the  Hamiltonian  algorithms
do not rely on being able to  sample from the marginal posterior distributions of the unknown parameters.
 Because of its relative ease in implementation, in this study we focus on the  MH  algorithm which we  introduce
 next.

To   numerically sample from the posterior
distribution in (\ref{postf}),
the Metropolis-Hastings MCMC algorithm proceeds as follows \cite{CS07, RCS14}.
Let  $\xi=\xi^i$ be a fixed parameter sample. Then, 
\begin{enumerate}
 \item[(1)] Generate a proposed sample $\xi^{\ast}$ from a proposal density 
 $q(\xi^{\ast}|\xi )$, and compute $q(\xi|\xi^{\ast} )$.
 \item[(2)] Set  $\xi^{i+1} = \xi^{\ast}$ with probability
 \be
 \label{accpr}
 \alpha(\xi^{\ast}|\xi) =
 \min\left[1, \frac{\pi(\xi^{\ast}|y_{obs})q(\xi^{\ast}|\xi ) }{\pi(\xi|y_{obs})q(\xi|\xi^{\ast} ) } \right].
 \e
Else, set $\xi^{i+1} = \xi$.
\end{enumerate}
Note that, at each step of the MH algorithm, the computation of  
$\alpha(\xi^{\ast}|\xi)$ in (\ref{accpr}) requires solving 
the discrete nonlinear system defined by
(\ref{disc1}) and (\ref{output}) to obtain $\pi(\xi^{\ast}|y_{obs}).$
 This happens especially when
 the proposal $\xi^{\ast}$ falls inside the parameter domain $\mathcal{D}$. Else,
$\pi(\xi^{\ast}|y_{obs}) = 0$,  and $\xi^{\ast}$ is immediately
discarded. In large-scale settings, the repeated solution of the  underlying discrete nonlinear PDE can
be quite computationally demanding. In the next section, we address 
this most expensive part of the MH scheme by deriving an efficient reduced-order model that will be used for the forward solve.

Now, denote by $\{\xi^i\}_{i=1}^{M}$ the samples generated by the MH algorithm.
These $M$  parameter samples
constitute a Markov chain.  Under mild regularity conditions, the Markov chain
can be shown to converge in distribution to the posterior density $\pi(\xi|y_{obs})$ \cite{RC2004}.
For a practical  implementation of the accept/reject rule in the MH algorithm, one typically computes a random
draw $\theta \sim \mathcal{U}(0,1),$ and then sets $\xi^{i+1} = \xi^{\ast}$ if 
$\theta < \alpha(\xi^{\ast}|\xi),$ and otherwise sets $\xi^{i+1} = \xi$.
For large-scale inverse problems, however, the computation of  $\alpha(\xi^{\ast}|\xi)$
is numerically unstable \cite{BSHL2014}. To circumvent this issue,  one can in this case  set $\xi^{i+1} = \xi^{\ast}$ if 
\[
 \ln{\theta}< \ln{ \pi(\xi^{\ast}|y_{obs}) } -\ln{ \pi(\xi|y_{obs}) }
           + \ln{ q(\xi^{\ast}|\xi )  } - \ln{q(\xi|\xi^{\ast} )  },
\]
or, otherwise, set $\xi^{i+1} = \xi.$

Next, note that the choice of the proposal density $q(\xi^{\ast}|\xi )$
is an extremely important part of the MH algorithm. We list here some  classes of proposal densities that are common in the literature.
One class of proposals is known
 as the independence proposals; they satisfy $ q(\xi^{\ast}|\xi^i ) = q(\xi^{\ast} )$.
An example of an independence proposal for MH is called the randomize-then-optimize (RTO) proposal
 which was introduced in \cite{BSHL2014}.
 The paper specifically discusses RTO in the context of nonlinear problems
  and shows that it  can be efficient compared to other proposals. The approach
 is based on obtaining candidate samples by repeatedly optimizing a randomly
perturbed cost function. 
Although the RTO yields relatively high acceptance rates for certain problems
 such as those considered in \cite{BSHL2014}, 
it is quite involved and very computationally expensive.

Another class of  proposals  consists of the  symmetric proposals;
that is,
\[
 q(\xi^{\ast}|\xi^i ) = q(\xi^i|\xi^{\ast} ), \quad 
 i= 1,2,\ldots, M.
\]
The consequence of this choice of $q(\xi^{\ast}|\xi )$  is that the acceptance 
probability becomes
\be
 \label{accpr2}
 \alpha(\xi^{\ast}|\xi^i ) =
 \min\left[1, \frac{\pi(\xi^{\ast}|y_{obs}) }{\pi(\xi^i |y_{obs}) } \right],
  \quad 
 i= 1,2,\ldots, M.
 \e
A typical example  is the so-called preconditioned Crank-Nicolson (pCN)
proposal density as discussed in \cite{CRSW13}.
The main challenge with the pCN is that the proposal density $q(\xi^{\ast}|\xi )$
typically needs to be tuned by the user to achieve an efficient MCMC method.
This requirement limits the usage of pCN in computational practice.

Another example of  symmetric proposal density, which we use in this work, is an adaptive Gaussian 
proposal introduced in \cite{HET2001}.
The proposal is defined recursively by
\be
\label{adapG}
q(\xi^{\ast}|\xi^{i} )  & \propto &
\exp\left( -\frac{1}{2}(\xi^{\ast} - \xi^{i-1})^T{\Gamma}_{i-1} (\xi^{\ast} - \xi^{i-1})  \right),
\e
where  ${\Gamma}_0=I,$ and 
\be
\label{covr}
{\Gamma}_i= \mbox{cov}\left( \{ \xi^j \}_{j=0}^{i-1}  \right) =\frac{1}{i}
\sum_{j=0}^{i-1} (\xi^j - \tilde{\xi})^T (\xi^j - \tilde{\xi}) + \epsilon I,
\e
with $\tilde{\xi} = i^{-1}\sum_{j=0}^{i-1}  \xi^j $ and $\epsilon$ a small positive
number, say, $\epsilon \approx 10^{-8}$. 
In the resulting  {\it adaptive Metropolis-Hastings} (AMH)  algorithm, we use lognormal proposal density,
in which case the covariance (\ref{covr})
is computed directly from the  log of the samples. Besides being  symmetric, this  proposal density 
enjoys relative ease in implementation. Hence,
in our numerical experiments, we implement AMH using the proposal density
(\ref{adapG}), together with the acceptance probability (\ref{accpr2}).

\section{Reduced-order modeling}
\label{rom_deim}
Reduced-order modeling involves deriving a relatively cheap and low-dimensional model
 from a typically high-fidelity and computationally costly model in such
  a way as to achieve  minimal loss of accuracy.
The aim here is  to reduce the cost of solving parameterized
discrete PDEs  at many parameter values. 
In 
reduced-order modeling,
the computation of the discrete problem (\ref{disc1}) proceeds in two steps: an
{\it offline} step and an {\it online}
step \cite{BMNP04, Cha2011, EF15, EF17}.
In this offline-online paradigm, the 
offline step constructs the low-dimensional approximation to the solution space,
and the online step uses this approximation -- the
{\it reduced basis} -- for the solution of a smaller reduced problem.
The resulting reduced problem provides an accurate estimate of the solution
of the original problem.
This section presents two reduced-order modeling techniques, POD-Galerkin and POD-DEIM-Galerkin approximations, 
for the general nonlinear system in (\ref{disc1}). The
key to the success of the reduced-order modeling techniques is for the cost of simulation
(that is, online operation count) to be independent of the
number of full-order degrees of freedom and the number of evaluations of the nonlinear term
$F(\cdot;\xi)$ in the full-order system.

\subsection{Constructing reduced-order model}
There are several techniques for constructing  reduced-order models \cite{BMNP04, Cha2011, EF15, EF17}. Perhaps
the commonest of these methods are the 
projection-based techniques.  They rely on Galerkin projection to construct a reduced-order system   that approximates
the original system from a subspace spanned by a reduced basis $Q $ of  dimension  $k\ll N. $
More precisely, 
assume that the columns of $Q = [\varphi_1,\varphi_2,\cdots,\varphi_k]\in \mathbb{R}^{N \times k}$ are the
basis vectors and the solution $u$ of the full model (\ref{disc1}) is adequately approximated by these vectors,
so that
\be
\label{qbasis}
u =Qu_r, 
\e
where
$u_r\in \mathbb{R}^k$ is the vector of reduced states. 
Then, using (\ref{qbasis}) and (\ref{disc1}) and applying Galerkin projection, one obtains the reduced model
\be
\label{disc2}
G_r(u_r;\xi) := A_ru_r + Q^TF(u_r;\xi) + B_r =0, 
\e
where $A_r = Q^TAQ\in \mathbb{R}^{k \times k}$ and $B_r=Q^TB\in \mathbb{R}^k.$
Similarly, from (\ref{output}), we get
\be
\label{output2}
y_r(u_r(\xi)) = C_ru_r(\xi), 
\e
where 
$C_r= CQ\in\mathbb{R}^{m\times k}$  and $y_r\in \mathbb{R}^m.$
In this work, we assume that the quantity of interest is the state
vector $u$ directly, so that $m=N $ and $C=I.$

The offline step is aimed at constructing
 the reduced basis $Q$, and it is the  most expensive part of the  computation.
Furthermore, the quality of the approximation is significantly affected by the  reduced basis $Q$.
There are several techniques for constructing the reduced basis; 
each of theses techniques is based on the fact that the solution lies in a 
low-dimensional space,
see e.g., \cite{Cha2011, GMNP07, HRS16, QMN2016}.
We use  the proper orthogonal decomposition (POD) method to construct the reduced basis.

Given a set of snapshots  of solutions obtained from the parameter space, the POD-Galerkin method constructs 
an optimal orthogonal (reduced) basis  in the sense  that  it minimizes the approximation error
with respect to the snapshots \cite{Rav05}.
More specifically, the POD method takes a set of $n_{trial}$
snapshots of the solution $S_U = [u(\xi^{(1)}), u(\xi^{(2)}),\cdots, u(\xi^{(n_{trial}) })]$
and computes the singular value
decomposition (SVD)
\be
\label{svd}
S_U =\bar{V}\Sigma W^T,
\e
where $\bar{V}=[\varphi_1,\varphi_2,\cdots,\varphi_{n_{trial}}]$ 
and $W$ are orthogonal and $\Sigma$ is a diagonal matrix with
the singular values sorted in order of decreasing magnitude. The reduced basis is
defined as $Q = [\varphi_1,\varphi_2,\cdots,\varphi_k]\in \mathbb{R}^{N \times k}$
with $k < n_{trial}$. This produces an orthogonal basis $Q$ which
contains the important components from the snapshot matrix $S_U.$ The major shortcoming
of POD is that it often requires an  expedient number of snapshots, $n_{trial}$,  to construct $S_U$.
It is possible that the number of solutions of the full model required to find a basis
with satisfactory accuracy could be quite large.

Despite the fact that the projection of the nonlinear operator $Q^TF(u_r;\xi)$ is of dimension $k$, 
it must be assembled at each step of a nonlinear iteration since it
depends on the solution. For a nonlinear solution method such as
the Newton method, each nonlinear iteration also
requires the construction of the Jacobian matrix associated with $F(u_r;\xi)$ and
multiplication by $Q^T$, where the costs of both operations depend on $N.$
More specifically, note that 
for the forward model defined in equation (\ref{disc1}), 
the Jacobian matrix  is
 \be
 \label{jacg}
 J_{G}(u;\xi) = A + J_{F}(u;\xi).
 \e
The Jacobian matrix $J_{G_r}(u)$ of the reduced model equation (\ref{disc2}) is
then
\be
\label{jacgr}
J_{G_r}(u_r;\xi) := Q^TJ_{G}(u_r;\xi)Q = A_r + Q^TJ_{F}(u_r;\xi)Q.
\e
This feature makes the POD-Galerkin  method 
inefficient in tackling the nonlinear problem.
This difficulty can be addressed  using the discrete interpolation
empirical method (DEIM), which we  discuss in Section \ref{deimsec}.

\subsection{The discrete empirical interpolation method (DEIM)}
\label{deimsec}
DEIM approximation is a discrete counterpart of the Empirical Interpolation Method (EIM) which was introduced in \cite{BMNP04}.
It is  used
 for  approximating
 nonlinear parametrized functions through sparse sampling.
Following \cite{BMNP04, Cha2011, EF17, Forstall2015},
we present next the POD-DEIM-Galerkin method.

Besides the reduced basis $Q\in\mathbb{R}^{N\times k}$ constructed from (\ref{disc1}), in the
DEIM framework,  a separate
basis is constructed to represent the nonlinear component of the solution. To construct this  basis,  one first computes
the snapshots of the nonlinear term
\[
S_F = [F(u(\xi^{(1)})), F(u(\xi^{(2)})),\cdots, F(u(\xi^{(s) }))], 
\]
 where $s \geq k$,
and then 
the SVD of the snapshot matrix:
\be
\label{svd2}
S_F = \tilde{V}\tilde{\Sigma} \tilde{W},
\e
where, as before,  $\tilde{V}=[v_1,\ldots,v_s]$  and $\tilde{W}$ are orthogonal, and  
 $\tilde{\Sigma}$ is a diagonal matrix  containing singular values. Then, 
using the first $n < s$ columns of $\tilde{V},$  one determines the DEIM basis
$V = [v_1,\ldots,v_n] \in \mathbb{R}^{N\times n}$.

Next, DEIM uses $n$  distinct interpolation points 
$\mathfrak{p}_1, \ldots , \mathfrak{p}_n \in 
\{1,\ldots , N\}$ and makes use of the DEIM interpolation  matrix
$P=[e_{ \mathfrak{p}_1 }, \ldots , e_{ \mathfrak{p}_n }] \in   \mathbb{R}^{N\times n}$
where $e_i$ 
is the $i$th  canonical unit vector with zeros in all components except $[e_i]_i =1$.
In \cite{CS10}, it is shown that the interpolation points and the interpolation  matrix
are constructed with a greedy approach using the DEIM basis $V$. 
To make the presentation self-contained, we show the construction in Algorithm \ref{deimalg}.

Observe  that  the $i$th interpolation point $\mathfrak{p}_i$ 
can be associated with the basis
vector in the $i$th column of the DEIM basis $V$. Moreover, 
the  DEIM interpolant of $F$ is defined by the tuple $(V,P)$, which  
is selected such that
the matrix
$P^TV\in\mathbb{R}^{n\times n}$  is nonsingular.
DEIM then approximates  the nonlinear function $F(u;\xi)$ by
\be
\label{deimap}
\bar{F}(u;\xi)  =  V(P^TV)^{-1}P^T F(u;\xi),
\e
where $P^T F(u;\xi)$ samples the nonlinear function at $n$ components only. 
This approximation satisfies $P^T\bar{F} = P^TF.$

Combining the DEIM and POD-Galerkin yields the POD-DEIM-Galerkin reduced
system
\be
\label{disc3}
G_{deim}(u_r;\xi) := A_ru_r + \bar{F}_r  + B_r =0, 
\e
where
\be
\label{disc3a}
\bar{F}_r(u_r;\xi) = Q^T \bar{F}(u;\xi) = Q^T V(P^TV)^{-1}P^TF(Qu_r;\xi). 
\e
The POD-DEIM-Galerkin reduced model approximation of the quantity of interest
 is defined  as in (\ref{output2}).
Moreover, the Jacobian $J_{G_{deim}}(u_r;\xi)$  of $G_{deim}(u_r;\xi)$  is
\be
\label{jacgdeim}
J_{G_{deim}}(u_r;\xi) :=  A_r + J_{\bar{F}_r}(u_r;\xi),
\e
where
\be
\label{jacfdeim}
J_{\bar{F}_r}(u_r;\xi) := Q^TJ_{\bar{F}}(u_r;\xi)Q =\boxed{Q^T V(P^TV)^{-1}}P^TJ_{F}(u;\xi)Q.
\e

\begin{algorithm}[t]
\caption{DEIM Algorithm \cite{CS10}}
\begin{algorithmic}[1]
\small
\State INPUT: $V = [v_1,\ldots,v_n] \in \mathbb{R}^{N\times n}$ linearly independent
\State OUTPUT: $\mathfrak{p} = [\mathfrak{p}_1, \ldots , \mathfrak{p}_n]^T\in \mathbb{R}^n,\; P \in   \mathbb{R}^{N\times n}$
\State $[|\rho|, \mathfrak{p}_1]=\max\{|v_1|\}$
\State $V=v_1, \; P=[e_{ \mathfrak{p}_1 }], \; \mathfrak{p} = [\mathfrak{p}_1] $
 \For{ $i=2$ {\bf to} $n$}
\State Solve  $(P^TV)c=P^Tv_i$  for $c$
\State $r = v_i - Vc$
\State $[|\rho|, \mathfrak{p}_i]=\max\{|r|\}$
\State $ V \leftarrow  [V \;\;  v_i], \; P\leftarrow [P \;\;  e_{ \mathfrak{p}_i }], \; \mathfrak{p} \leftarrow 
\left[\begin{array}{c}
\mathfrak{p} \\
\mathfrak{p}_i \\
\end{array}\right].
$
\EndFor
\end{algorithmic}
\label{deimalg}
\end{algorithm}

Solving (\ref{disc3}) with  Newton's 
method only requires the evaluation of entries of  the nonlinear function $F$ 
at indices associated with the interpolation points given by
$P$, instead of at all $N$ components. 
Note that, for the costs to be low, for each $i$, $F_i$ should depend on $\mathcal{O}(1)$
values of the solution $u,$ and similar requirements are needed for the Jacobian. 
This condition is valid for a typical (finite difference or finite element) discretization
of PDEs \cite{AHS13}.
The corresponding computational procedure of
the DEIM method is split into an offline phase where the DEIM
reduced system is constructed and an online phase where it is evaluated. 
The offline phase obviously involves high computational costs; however, these costs are incurred only once and
are amortized over the online phase.
In particular, the construction of  the matrix $Q^T V(P^TV)^{-1}$ highlighted in  (\ref{jacfdeim})   is part of the offline
phase.

In what follows, we refer to the combination of AMH and the DEIM model as the AMH-DEIM scheme. Similarly,
we write AMH-Full scheme for the combination of AMH and the full model.

\section{Numerical experiments}
\label{numexp}
In this section, we present  results on the performance of  AMH-DEIM  and AMH-Full schemes for a statistical
inverse problem, as well as 
the reduced-order  (POD-Galerkin and DEIM) models for a nonlinear forward problem. To this end, consider the following nonlinear
 diffusion-reaction problem posed in a two-dimensional spatial domain \cite{Cha2011, GMNP07}
\begin{eqnarray}
\label{nonlinear}
-\nabla^2 u(x_1,x_2) + F(u(x_1,x_2) ; \xi) &=& 100 \sin(2\pi x_1) \sin(2\pi x_2),\\
 F(u; \xi)  &=& \frac{\xi_2}{\xi_1}\left[\exp(\xi_1 u) - 1\right],\nonumber
 \end{eqnarray}
where the spatial variables $(x_1, x_2) \in \Omega = (0, 1)^2$ and
the parameters are $\xi = (\xi_1, \xi_2)\in\mathcal{D} = [0.01, 10]^2 \subset \mathbb{R}^2$,
with a homogeneous Dirichlet boundary condition.

Equation   (\ref{nonlinear}) is  discretized on a uniform mesh in $\Omega$ with $32, \;64$  and  $128$ 
grid points in each direction using
centered differences resulting, respectively, 
in  $N=1024,\; 4096$ and $16384$ (unless otherwise stated). We solved the resulting 
system of nonlinear equations  with inexact Newton-GMRES method 
as described in \cite{CTK95}.
Recall that, for a fixed $\xi,$  an inexact Newton method solves for $u$ in the nonlinear problem (\ref{disc1})
by approximating the vector
  $z$ in the equation for the Newton step
\begin{eqnarray}
\label{newtdir}
J_G(u^{(i)};\xi)z= (A + J_{F}(u^{(i)};\xi))z = -G(u^{(i)};\xi),
 \end{eqnarray}
 with
 \[
  u^{(i+1)}=u^{(i)} + z, \;\;
 i=0,1,2,\ldots,
 \]
such that
\begin{equation}
\label{newtcond}
 || J_G(u^{(i)};\xi)z + G(u^{(i)};\xi) ||\leq \eta_i ||G(u^{(i)};\xi)||,\;\;
 i=0,1,2,\ldots,
\end{equation}
where the parameters $\{\eta_i\}$ are the so-called forcing terms. To realize
the inexact Newton condition (\ref{newtcond}), Newton iterative methods typically
 apply an iterative method (GMRES method in this work \cite{book::Saad}) to the equation for the Newton step and 
terminate
that iteration when (\ref{newtcond}) holds. We  refer to this linear
iteration as
the {\it inner iteration} and,  the nonlinear iteration 
as the {\it outer iteration}.

If the forcing terms $\{\eta_i\}$ 
in the inexact Newton method are uniformly (and strictly) less than $1$,
then the method is locally convergent \cite{DES82}. 
For more details of local convergence theory and the
role played by the forcing terms in inexact Newton methods, see e.g., \cite{AML2007, EW96, CTK95}. 

The  forcing terms  are  chosen in such a way as to solve the linear equation for
the Newton step to just enough precision to make good progress when far from a
solution, but also to obtain quadratic convergence when near a solution \cite{CTK95}. 
Following \cite{CTK95}, in our computations we  chose 
\begin{equation}
\label{eta}
 \eta_i = \min(\eta_{max},\max(\eta_i^{safe},0.5\tau/||G(u^{(i)};\xi)||)),
\end{equation}
where 
the parameter $\eta_{max}$ is an upper limit on the forcing
term and
\begin{align}
\label{eta2}
 \eta_i^{safe} &=
  \begin{cases}
    \eta_{max}        & i =0, \\
   \min(\eta_{max},\eta_i^{res} )        & i>0,\; \gamma\eta^2_{i-1}\leq 0.1,\\
   \min(\eta_{max}, \max(\eta_i^{res},\gamma\eta^2_{i-1} ) )        & i>0,\; \gamma\eta^2_{i-1} > 0.1,
  \end{cases}
\end{align}
with
\[
\eta_i^{res} = \gamma || G(u^{(i)};\xi)||_2/|| G(u^{(i-1)};\xi)||_2,
\]
the constant $0.1$ being somewhat arbitrary, $\gamma \in [0,1)$ and $\tau = \tau_a + \tau_b||G(u^{(0)};\xi)||.$
In our experiments, we set
$\gamma =0.9, \;\eta_{max}=0.25,$
$ \tau_a = \tau_b =10^{-6},$
and
 $u^{(0)} \equiv 0.$  
We use the same strategy for  the DEIM reduced model 
with $G$ in (\ref{disc1}) replaced by $G_{deim}$ in (\ref{disc3}).

We compare the performance of the solvers by plotting the relative
nonlinear residuals from the full and DEIM models. More precisely, we plot
\begin{itemize}
 \item  the relative nonlinear residual  $R_{full}:=|| G(u^{(i)};\xi)||_2/|| G(u^{(0)};\xi)||_2$
 against the number of  calls of the function $G$ in the  full model defined by (\ref{disc1}).
 \item   the relative nonlinear residual  $R_{deim}:=|| G_{deim}(u_r^{(i)};\xi)||_2/|| G_{deim}(u_r^{(0)};\xi)||_2$
 against the number of calls of the function  $G_{deim}$, in  the DEIM model given by  (\ref{disc3}).
\end{itemize}
In particular, we make $R_{full} < \tau$
and $R_{deim} < \tau_r,$ where we choose $\tau_r$ small enough with $\tau_r < \tau$,
such that we expect the approximation $Qu_r$ to be as good as the numerical solution $u.$
In our numerical experiments, we specifically set $\tau_r = 10^{-2}\tau.$

\subsection{Preconditioning}
Now, as we have mentioned earlier, at each Newton step in both the full and reduced models, we elect to solve  the associated linear system 
with the GMRES method. However, Krylov solvers (including GMRES method) are generally slow and require appropriate
preconditioners to reduce the associated computational cost \cite{ESW14}.
For an arbitrary linear system of equations  $AX = B,$  we recall here that a preconditioner 
is, loosely speaking, a matrix $\mathcal{M}$ such that
(a) the system  $\mathcal{M}A X = \mathcal{M}B$ yields a matrix $\mathcal{M}A $
whose eigenvalues are well clustered, and 
(b) the matrix-vector product $\mathcal{M}Z,$ is relatively easy to compute
for any vector $Z$ of appropriate dimension. In the context of Krylov solvers,
the consequence of (a) is that the solver converges quickly 
\cite{ESW14}. In particular, to solve (\ref{newtdir}),  one could precondition 
the linear system with an approximation $\mathcal{M}$ to
the action of $A^{-1}$ (that is, an  approximate Poisson solver).

In this work, we use the AGgregation-based algebraic MultiGrid iterative method
 (AGMG) as a preconditioner
 \cite{NN2012,  NY2012, NY2010}.
AGMG employs  algebraic multigrid method to solve  the linear system of equations  $AX = B$.       
 The algebraic multigrid method is  accelerated either by the Conjugate Gradient method (CG) or by
        the Generalized Conjugate Residual method (GCR, a variant of GMRES). 
AGMG is implemented as a ``black box''. In particular, since $A$ (the Laplacian, in our case) is symmetric and positive
definite, in our numerical experiments,
we employ the CG iterations option in the software\footnote{For our problem,
there were $10$ CG iterations, and the CG method stopped when the Euclidean norm of the 
residual relative to that of the right hand side was less than $10^{-6.}$};  it uses one
V-cycle of AMG with symmetric Gauss-Seidel (SGS) pre- and post-smoothing sweeps to
approximate the action of $A^{-1}.$
We apply the preconditioner  to the full model (\ref{disc1}) as follows:
\be
\label{disc1prec}
u + \mathcal{M}F(u;\xi) + \mathcal{M}B =0. 
\e
Note here that the preconditioner $\mathcal{M}$ in (\ref{disc1prec}) is used as a Poisson
solver.
For semilinear boundary value problems\footnote{That is,
our model is linear in the highest
order derivative, which is the Laplacian.}, preconditioning  with an exact inverse for
the high-order term typically
yields optimal convergence
of the linear iteration in the sense that 
the iteration counts are independent of the mesh spacing \cite{MP90}.

It is a common computational practice to solve the reduced problem  using direct methods 
\cite{BMNP04, Cha2011}. 
However, as  observed in \cite{EF15, EF17}, the size of the reduced
system corresponding to a prescribed accuracy is essentially determined by the dynamics of the problem.
Quite often,  although the size of the reduced model is significantly smaller than that of the full model,
solving the former with direct solvers can still be computationally expensive.
In this case, the use of
preconditioned iterative
methods may be more effective to solve the reduced problem. 
To this end, we precondition the DEIM model  (\ref{disc3}) as follows:
\be
\label{disc3prec}
 u_r + \boxed{ Q^T \mathcal{M} V(P^TV)^{-1}}P^TF(Qu_r;\xi)  + \boxed{ Q^T \mathcal{M}B} =0, 
\e
where the  matrices $Q^T \mathcal{M} V(P^TV)^{-1}$  and $Q^T\mathcal{M} B$
should be  precomputed.  The  preconditioned DEIM model (\ref{disc3prec}) implies that 
the action of $A_r^{-1}=(Q^TAQ)^{-1}$ is 
approximated with $Q^T\mathcal{M}Q$. 
This preconditioning strategy is  inspired by the experimental results
in \cite{EF15}, which show that all the eigenvalues of the matrix
\[
(Q^TA^{-1}Q)A_r=Q^TA^{-1}QQ^TAQ
\]
 are bounded below by $1$, and the largest eigenvalues grow
only slightly with spatial dimension. 
This finding suggests that the condition number of the preconditioned reduced matrix is only weakly
dependent on the spatial mesh.
Note, however, that in general,
$(Q^TAQ)^{-1} \neq Q^TA^{-1}Q$ since the matrix $Q$ is often tall and
rectangular\footnote{An alternative 
way to implement  the action of $A_r^{-1}$ is by 
computing the $LU$ decomposition of $Q^TAQ$, and then 
using the factors of the inverse of the resulting
matrix \cite{EF17}; that is,
\[
 A_r^{-1}=(Q^TAQ)^{-1} = (LU)^{-1}= U^{-1}L^{-1}.
\]
}.

 Next, observe that  since (\ref{nonlinear}) is  discretized  on a uniform $5$-point 
centered difference mesh, each row of the  discrete Laplacian has at most 
$5$ nonzero elements. Thus,  the computational cost of matrix-vector product 
by  $\mathcal{M}$ required  an iterative solver in the  preconditioned full model (\ref{disc1prec}) is 
$\mathcal{O}(N).$
On the other hand,  the cost of analogous computation with 
$Q^T \mathcal{M} V(P^TV)^{-1}P^T$  in the preconditioned DEIM model (\ref{disc3prec})
is $\mathcal{O}(k^2).$ 
Hence,  if $k \ll N,$  the computational complexity of the forward model (\ref{nonlinear}) 
is significantly reduced when  (\ref{disc3prec}) is used in place of   (\ref{disc1prec}).

\subsection{Solving the parameterized nonlinear system}

\begin{figure}[bt]
\centering
\includegraphics[width=1.05\textwidth,height=0.55\textwidth]{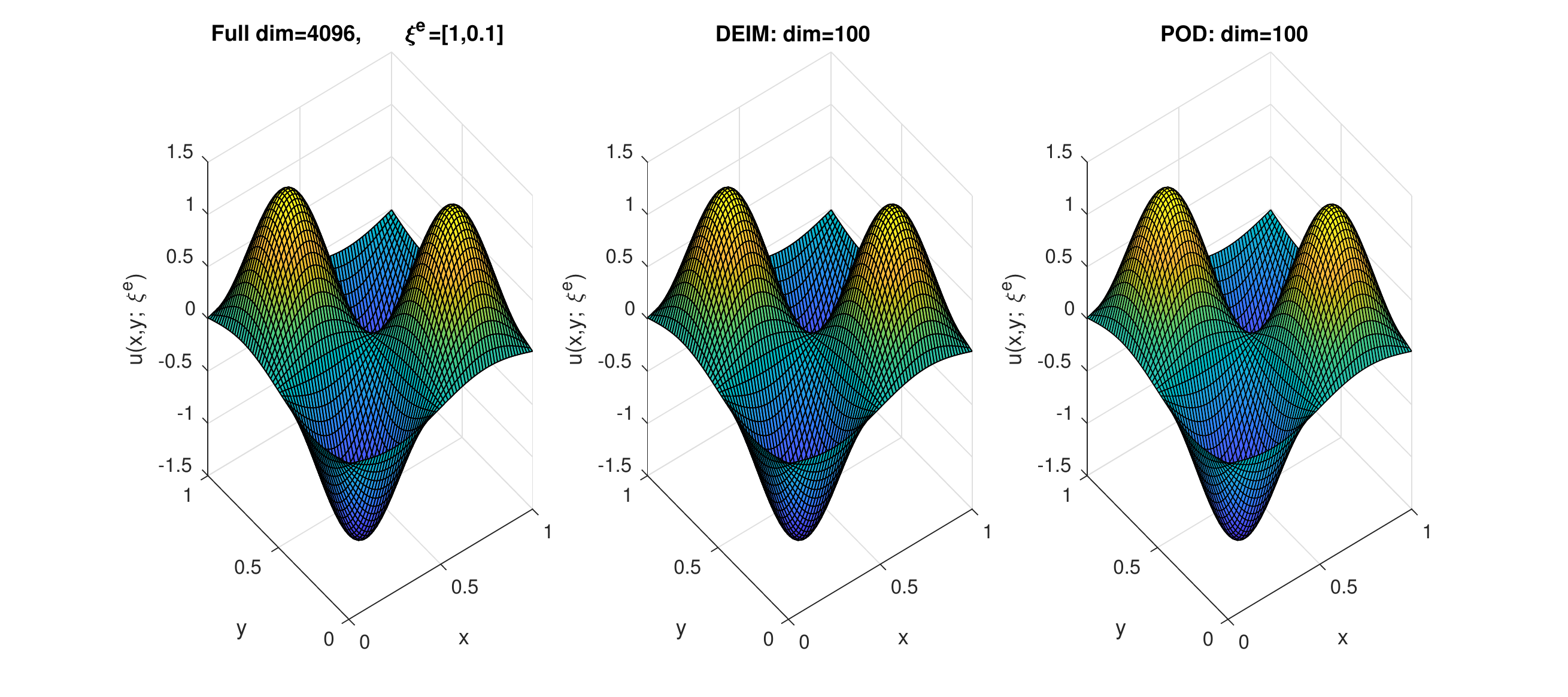}
\caption{Numerical solution   at $\xi  = [1, 0.1]$ from the full-order forward model (left, $ N=4096$) 
together with  the solutions from
 reduced systems: DEIM (middle) and POD-Galerkin (right).
}
\label{forward_sols101}
\end{figure}

\begin{figure}[bt]
\begin{minipage}{.51\linewidth}
\centering
\subfloat[]{\label{main:a}\includegraphics[scale=.475]{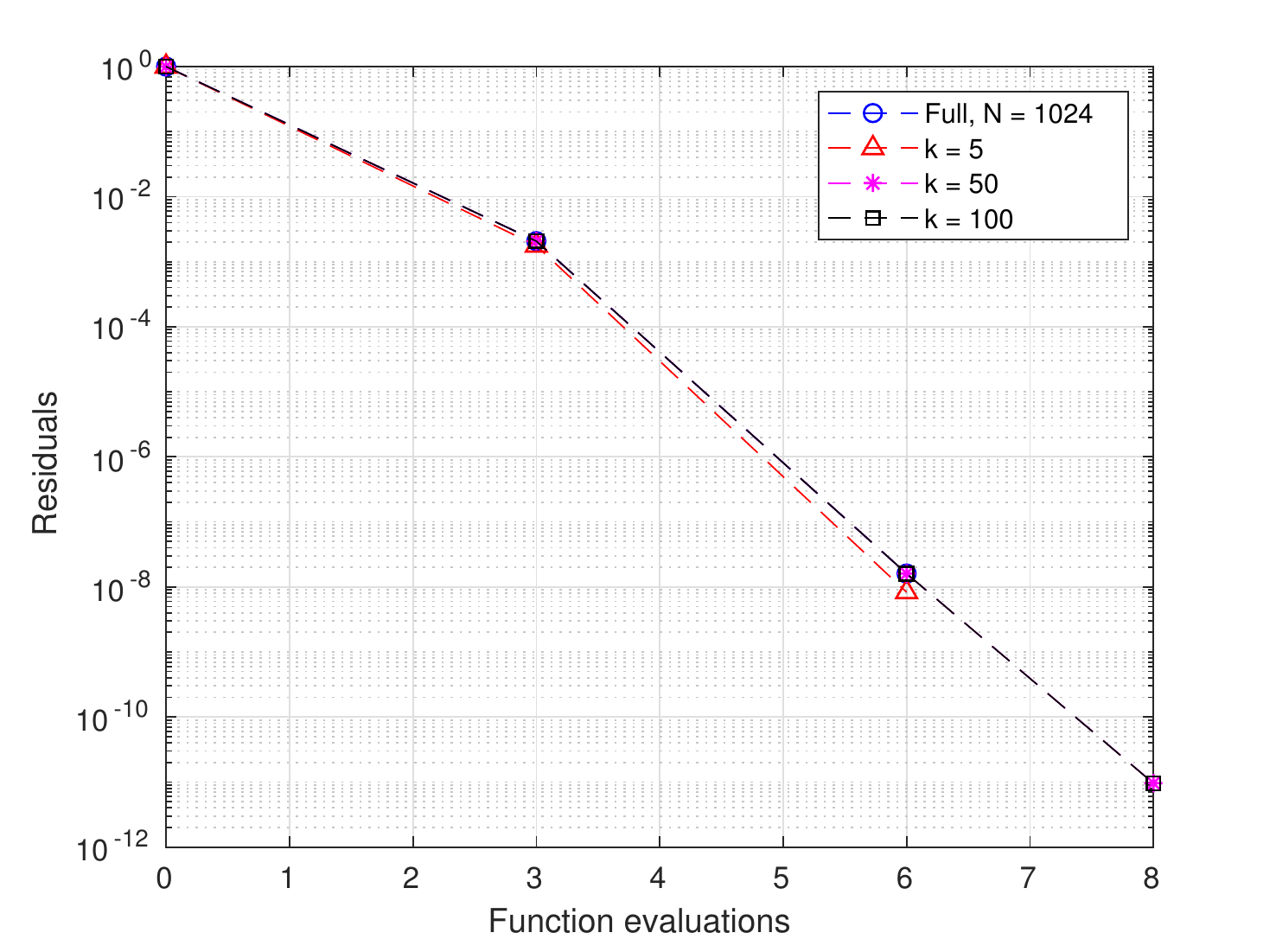}}
\end{minipage}%
\begin{minipage}{.51\linewidth}
\centering
\subfloat[]{\label{main:b}\includegraphics[scale=.475]{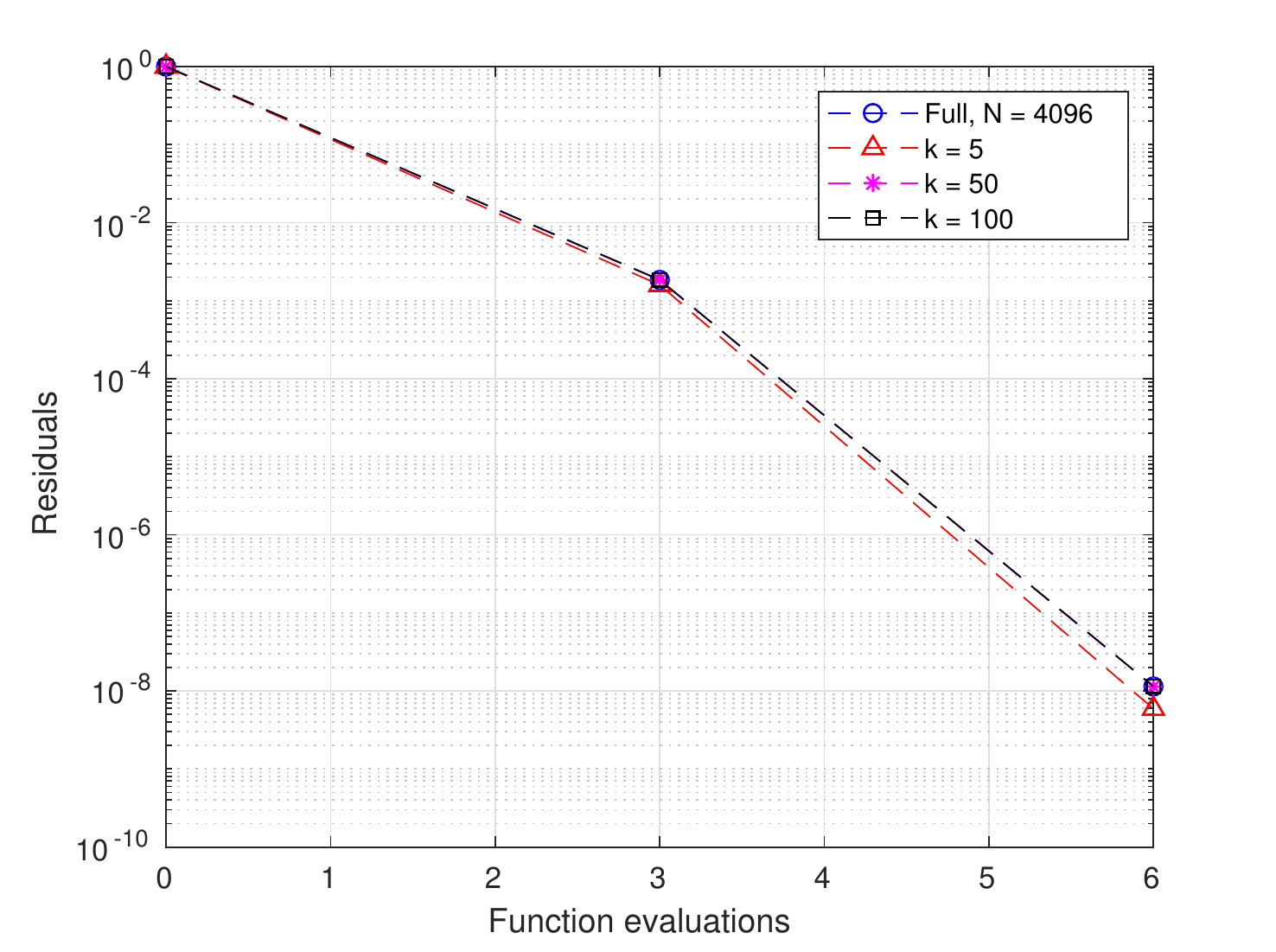}}
\end{minipage}\par\medskip
\centering
\subfloat[]{\label{main:c}\includegraphics[scale=.475]{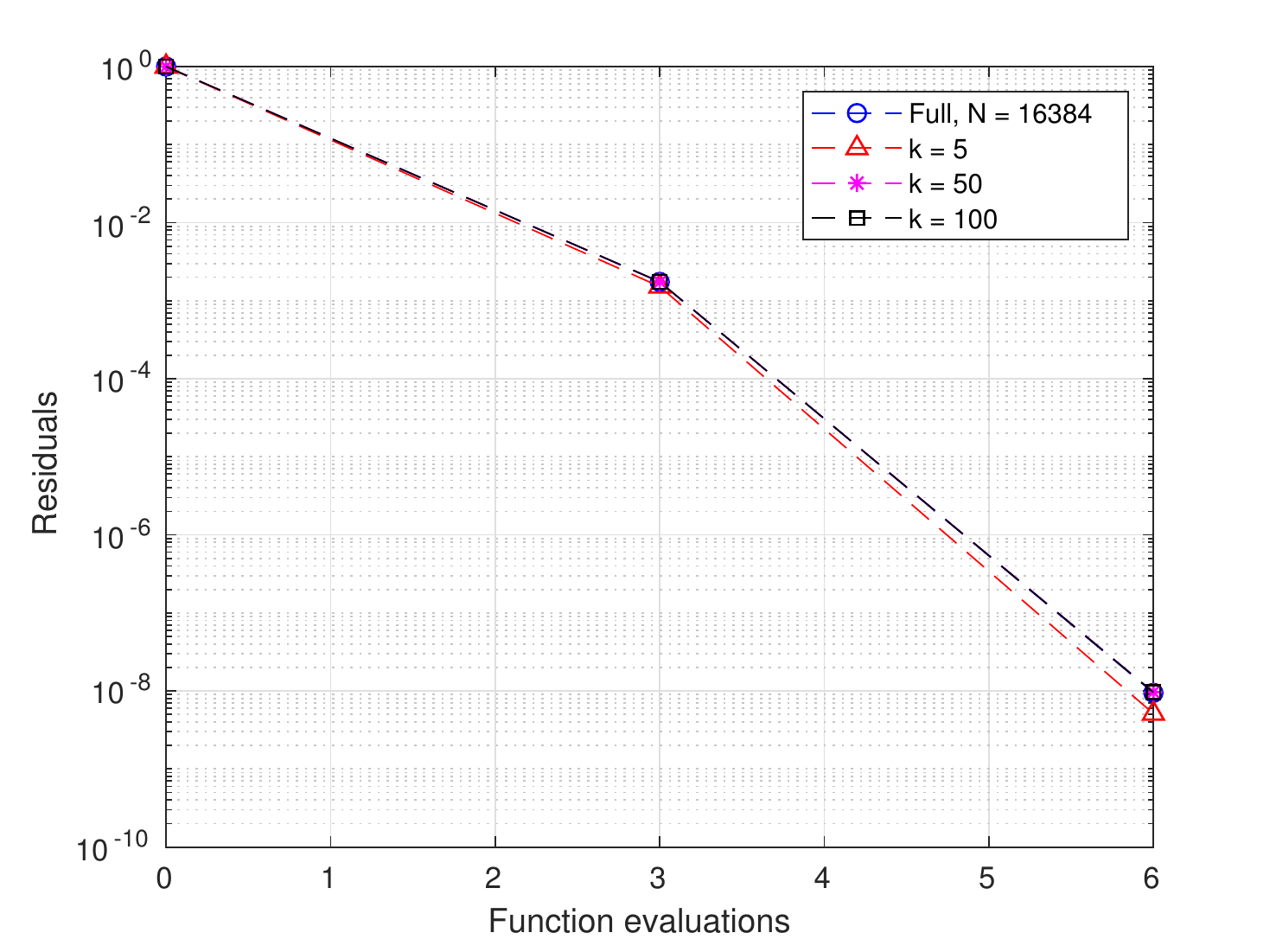}}
\caption{Performance of the full-order  model for $(a)\; N=1024, \; (b)\; N=4096, \; (c)\; N=16384,$
and  DEIM models for $k=5,\;  50,\; 100$ (for a fixed $N$) with the forward problem 
at $\xi  = [1, 0.1]$. Counts of function evaluations (fe) corresponding to
nonlinear iterations are indicated by circles for the full model and triangles, square and stars for the DEIM models.
}
\label{forw_nk325}
\end{figure}

To solve the Bayesian statistical inverse problem, we first discuss how we evaluated the forward model. We
 generated a reduced basis with $Q\in \mathbb{R}^{N\times k}$
and DEIM basis $V\in \mathbb{R}^{N\times n}$. In all our experiments, we set $k=n$.
We used $100\times 100$ grid points in  $\mathcal{D}$ which correspond to $n_{trial} =10000$ snapshots.

Additionally, we use the random sampling 
technique\footnote{Another technique for constructing the reduced bases for 
DEIM is the greedy algorithm \cite{BMNP04,Boyaval10}.
This method is significantly more
time-consuming than the random sampling approach for the problem
considered herein.}
as described 
in \cite{EF17, Forstall2015} to 
construct $Q$ and $V.$
This approach proceeds with a random
sample of $n_{trial}$ parameters, and  a snapshot is taken only if the reduced solution
at the current sample fails some error criterion. 
The basis is initialized using a single snapshot. Then,
for each of the $n_{trial}$ parameters, the reduced problem
is solved. If the error indicator $\eta_{\xi}$ of the reduced solution is below a prescribed  tolerance $\tau_d$, the
computation proceeds to the next parameter. Otherwise, the full model is solved
and the snapshot is used to augment the reduced basis. In 
our  experiments, we set $\tau_d= 10^{-4}$
and, measure the quality of the reduced model using the full residual since it is an easily computable quantity that indicates
how well the reduced model approximates the full solution. More specifically,
we use  the 
error indicator 
\[
\eta_{\xi} =   || G(u;\xi)||_2/||B||_2.               
 \]

We generated a data set by using
 (\ref{obs}) with  the ``true parameter'' $\xi =[1,0.1]$ and 
a Gaussian noise  vector $\varepsilon$ with $\sigma =  10^{-2}.$
First,  Figure  \ref{forward_sols101},  shows the noise-free solution
of the forward problem at the ``true parameter'' $\xi =[1,0.1]$ using the full model, as well as 
those solutions obtained with the  reduced models at $\xi$. We have also computed the following relative errors:
\begin{eqnarray}
 rel_{deim}: &=& \frac{||u - Qu_{deim}||_2}{||u||_2} = 3.2603\times 10^{-6},\nonumber  \\
 rel_{redb}: &=& \frac{||u - Qu_{redb}||_2}{||u||_2} = 1.9851\times 10^{-7},\nonumber \\
  rel_{dr}: &=&  \frac{||u_{deim} - u_{redb}||_2}{||u_{redb} ||_2} =  3.2543\times 10^{-6},\nonumber
\end{eqnarray}
where $rel_{deim}$ denotes the relative error with respect to the full and DEIM solutions, 
$rel_{redb}$ denotes the relative error with respect to the full and POD-Galerkin solutions, and 
$ rel_{dr}$ denotes the relative error with respect to the DEIM and POD-Galerkin  solutions. Observe here
that the POD-Galerkin  solution yields slightly more accurate solution than the DEIM solution. However,
since, as we have noted earlier, the computational cost of solving the forward problem using the POD-Galerkin  reduced model grows
with $N,$ we  use the DEIM reduced model in the rest of our computations.

Next, Figure \ref{forw_nk325} depicts the performance of 
the preconditioned full model compared to that of the preconditioned  DEIM model. In  particular, 
the size of the full model $N$ in each of the figures is fixed,
while the reduced basis dimension $k$ and DEIM basis dimension $n,$ are varied with $n=k;$ that is,  
$k=n=5, 50, 100.$ Plots $(a)$, $(b)$ and  $(c)$ correspond to
$N=1024, \; 4096$ and $16384,$ respectively.  They show
the relative nonlinear residual  $|| G(u^{(i)};\xi)||_2/|| G(u^{(0)};\xi)||_2$  plotted
 against the number of function $G$  evaluations required by all inner and nonlinear
iterations to compute $u^{(i)}$ for the  full model in (\ref{disc1}),  
together with
the relative nonlinear residual  $|| G_{deim}(u_r^{(i)};\xi)||_2/|| G_{deim}(u_r^{(0)};\xi)||_2$
 plotted against the number of function $G_{deim}$ evaluations required by all inner
 and nonlinear
iterations to compute  $u_r^{(i)}$ for  the DEIM model in (\ref{disc3}).

 \begin{figure}[bt]
\centering
\includegraphics[width=0.85\textwidth,height=0.50\textwidth]{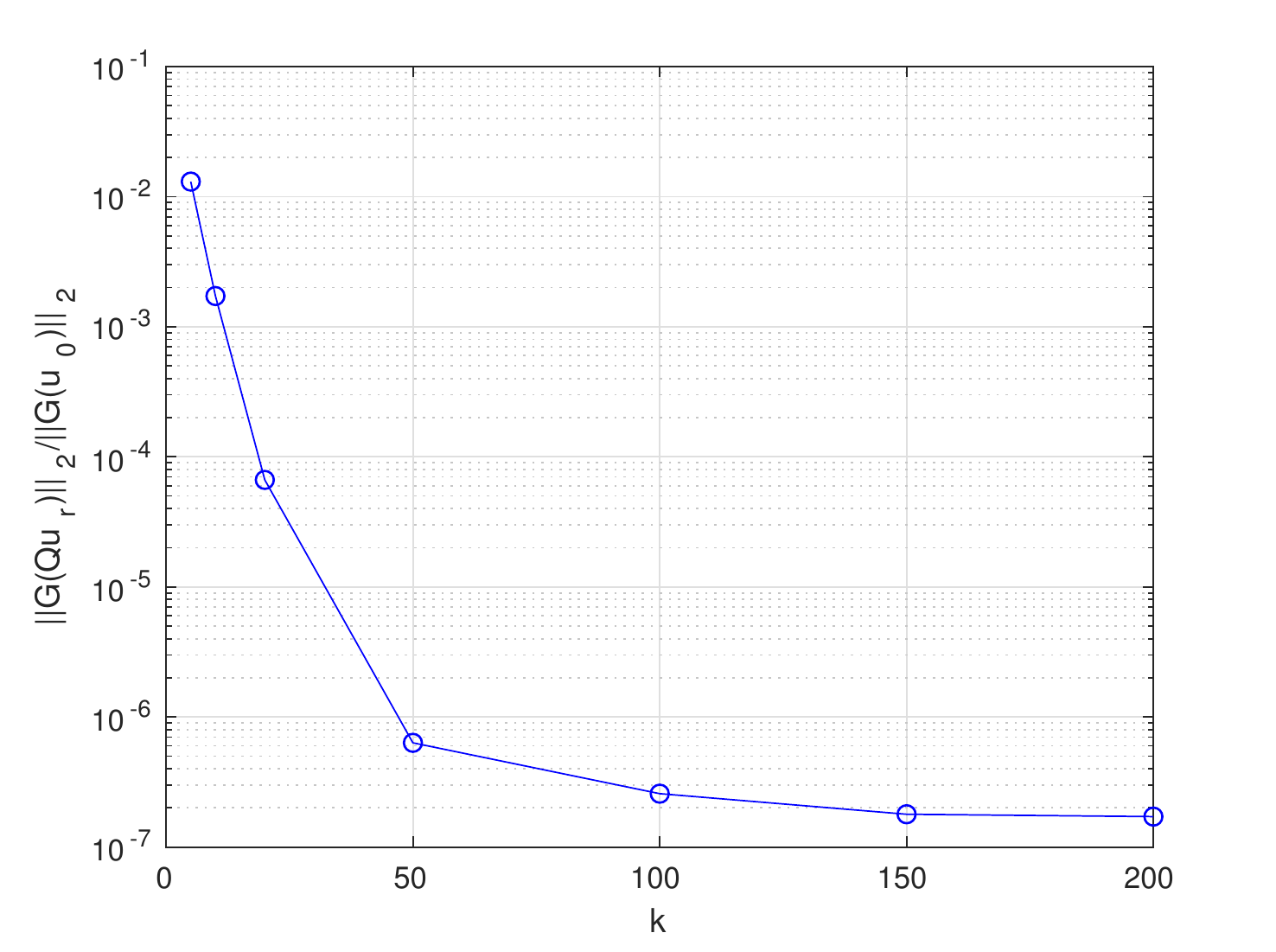}
\caption{Impact of the size of DEIM and reduced bases $k$ on the accuracy of the DEIM model. Results
were obtained with $N=16384.$
}
\label{deimac}
\end{figure}

The counts of function evaluations corresponding to
nonlinear iterations in these plots are indicated by small circles 
(for the full models) and triangles, square, star (for the DEIM models).
With the exception of the plot for $k=5,$ the plots for the DEIM models 
(i.e., $k=50$ and $k=100$)
coincide with the plot for all the full models for each  $N.$
From these plots,    one can compare both
 the number of nonlinear iterations and  the total cost. 
Observe from  these plots, that for a fixed $N,$  as 
the size of the DEIM model $k$  increases, 
the number of nonlinear iterations  and function evaluations
 required to solve  both the preconditioned full model and DEIM models
remain fairly constant; that is, about $3$ to $4$ nonlinear iterations
and $6$ to $8$ fe. 

To test the accuracy of the model reduction techniques
for the forward problem, reduced basis and DEIM basis of various sizes, $k$
and $n$, respectively, (with $k=n$) were constructed using  $n_{trial} =625$ snapshots. 
To do this, we computed and plotted the relative residuals
\[
 rels:= \max_{\xi}|| G(Qu_r;\xi)||_2/|| G(u_0;\xi)||_2, \;\;\; \forall \xi \in \mathcal{E},
\]
where $\mathcal{E}$ is the set of $n_{trial} $ parameters. In Figure \ref{deimac},
 $rels$ are plotted against $k$.
The figure shows that after $k=150$, the accuracy of the DEIM model does not improve by increasing 
$k.$

Table \ref{relerr_time}, in combination with Figure \ref{deimac}, shows
that the computational time for solving 
the  DEIM models is generally much smaller 
than that for solving the full model. In particular, although the 
cost of DEIM 
increases with the number of vectors $k,$  similar accuracy is obtained by DEIM
 for $k\geq 50$ at a significantly   lower cost.
\vspace{2mm} 
 \begin{center}
\begin{tabular}{l|llllll}
\hline
$ N$    &  $k=5$         &  $k=20$           &  $k=50$        & $k=100$       &  full        &  \\
\hline
\hline
$1024$  &  $(0.042, 6)$  &   $(0.068, 8)$    & $(0.082, 8)$   & $(0.140, 8)$  &  $(0.056, 6)$ & \\
\hline
$4096$  &  $(0.060, 6)$  &   $(0.073, 6)$    & $(0.087, 6)$   & $(0.154, 6)$  &  $(2.230, 6)$ & \\
\hline 
$16384$ &  $(0.057, 6)$  &   $(0.069, 6)$    & $(0.079, 6)$   & $(0.142, 6)$  &  $(10.45, 6)$ & \\
\hline
\end{tabular}
\captionof{table}{CPU times   and  the number of functional evaluations (t, fe)
needed 
 to compute the solution at $\xi  = [1, 0.1]$ by the preconditioned models
with $k (=n)$ and $N$ representing the dimensions of the DEIM  and full models respectively.
}
\label{relerr_time}
\end{center}
For example, in Table \ref{relerr_time},  the computational times required by the 
full model increased  from $0.056s$ for $N=1024$ to  $10.45s$ for $N=16384$ 
while the DEIM model  required a maximum of $0.154s$ with benign dependence on $k,$
but, as expected, its costs were independent of $N.$
To solve the  statistical inverse problem, we set  $k=n=100$ 
in all   subsequent computations with the AMH-DEIM model in Section \ref{statinv}.

 \begin{remark}
Several attempts have been made in recent years  to improve the performance of standard DEIM as presented 
in this work; see e.g.,  \cite{DG16, PBWB14, PDG18, PW15,  KS19} and the references therein.
A major strategy in this direction is the so-called {\it randomization}  technique \cite{DG16, PDG18,  KS19}.
 Originally introduced in \cite{DG16}, the randomized
sampling approaches for computing  DEIM indices have been analyzed in \cite{PDG18,  KS19}.
In particular, the paper \cite{PDG18} proposes {\it overdetermined} DEIM (ODEIM), which employs
randomized oversampling to tackle stability issues associated with the use of standard DEIM. We implemented ODEIM 
and found that although it slightly improves the accuracy of the reduced model, the algorithm yielded increased
 computational times  for  solving our forward problem compared to those from standard DEIM as reported in Table \ref{relerr_time}.
 \end{remark}

\subsection{Solution of the statistical inverse problem} 
\label{statinv}
Recall that the goal of solving the inverse problem
is to numerically approximate the posterior probability distributions of the  parameters of interest
$\xi_1$  and $\xi_2$, as well as to  quantify uncertainty  from
the distribution. To achieve this, we compute $M=20000$ steps of AMH 
as described in Section \ref{mcmc}. Here, we use an initial Gaussian 
proposal with covariance ${\Gamma}_0=I$ and update the covariance using 
(\ref{covr}) after every $100$th sample. After discarding the first 
$10000$ samples due to {\it burn-in}\footnote{This is the initial
stage of any MCMC chain, when the elements of the chains stagger
from their starting values to the region of relatively
high probability of the target density.
It is important to discard the  
samples corresponding to burn-in to avoid computing biased statistics from the 
MCMC chain. The remaining MCMC chain after discarding  the burn-in samples
is said to be in equilibrium,
 and provided it is long enough, it can be treated as a collection
 of samples from a target density \cite{BSHL2014}.}, we plot the results 
 from AMH-DEIM model (with $k=n=100$) and AMH-full model (with $N=16384$).

 We assess the performance of AMH-DEIM model and AMH-full model in
sampling the target density via the following criteria:
\begin{itemize}
\item Autocorrelation: 
The efficiency
with which the MCMC method explores the target density is quite important. It is determined by
the degree of correlation in the 
$(\xi_1, \xi_2)$ AMH chains. More specifically,
suppose $\{\delta_j\}_{j=1}^J$  is a Markov chain generated by the AMH algorithm,
where  $\{\delta_j\}$ are identically distributed with variance $\sigma^2.$
 Suppose also that the covariance is translation invariant, that is,
 $\mbox{cov}(\delta_j,\delta_{k})=\mbox{cov}(\delta_1,\delta_{1+|j-k|}),$
 for $1\leq j,k \leq J,$  where  $\mbox{cov}(\cdot,\cdot)$ denotes the covariance function.
Then, the autocorrelation function (ACF) $\rho$ of the  $\delta$-chain is given by
\[
 \rho_j = \mbox{cov}(\delta_1,\delta_{1+|j|})/\sigma^2. 
\]
 The 
ACF decays very fast to zero as $J\rightarrow \infty.$
Moreover, the integrated autocorrelation time, $\tau_{int}$, (IACT)
of the chain is defined as
\begin{eqnarray}
 \tau_{int}(\delta)&:=&\sum_{j=-J+1}^{J+1}\rho_j\nonumber\\
 &=&\sum_{j=-J+1}^{J+1}\mbox{cov}(\delta_1,\delta_{1+|j|})/\sigma^2\nonumber\\
 &\approx& 1 + 2\sum_{j=1}^{J-1}\left(1-\frac{j}{J}\right)\mbox{cov}(\delta_1,\delta_{1+j})/\sigma^2,
 \end{eqnarray}
where, in practical  implementation, $J$ is often taken 
to be $\left\lfloor{10\log10(M)}\right\rfloor,$ \cite{BSHL2014}.
If $\mbox{IACT} = K,$ say, then it means that the roughly every $K$th sample 
of the $\delta$ chain is independent.
\item Geweke test: This uses the Central Limit Theorem to compute 
a statistic $R_{Geweke}$ which is used to obtain the probability $p$
of accepting the null hypothesis of, say, $\delta_{10\%}=\delta_{50\%},$
where $\delta_{10\%}$ and $\delta_{50\%},$  denote the means of
the first $10\%$ and last $50\%$ of the chain, respectively. It should be noted that
a high value of $p$ suggests that the chain is in equilibrium.
\item Confidence interval (CI): This shows the interval within which all the 
samples of the chain lie with a prescribed probability. The midpoint of this interval
(which is the mean of these
samples) can be used to estimate  the true parameter.
\end{itemize}

 \begin{figure}[bt]
\centering
\includegraphics[width=1.0\textwidth,height=0.650\textwidth]{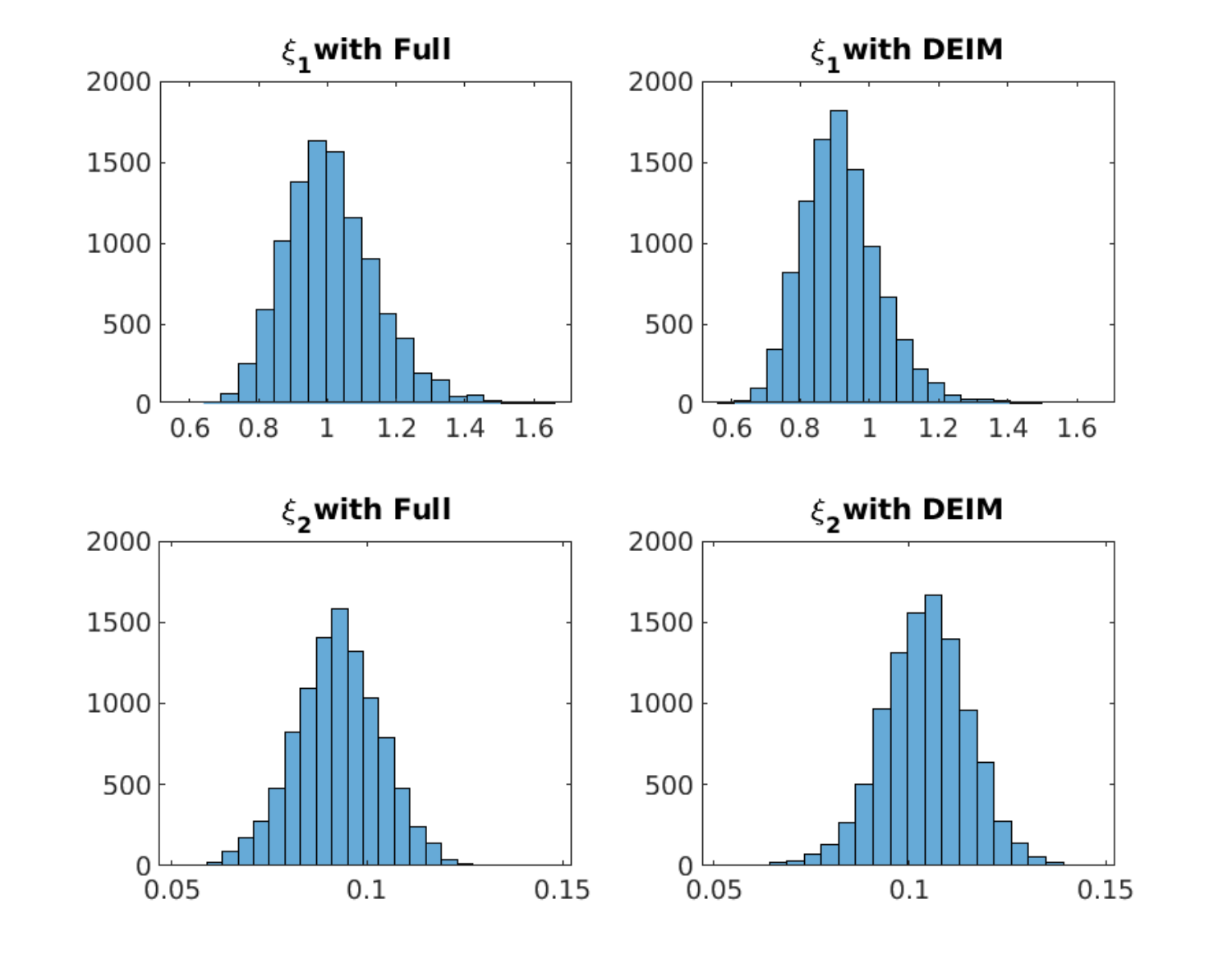}
\caption{Histograms of the  posterior distributions  for the parameters $\xi  = [\xi_1, \xi_2].$
They were obtained from AMH-Full (left) and  AMH-DEIM (right) models with  $M=20000$ MCMC samples.
 With AMH-DEIM model, the $95\%$ CI is $ (0.75923, 1.18712)$ for $\xi_1$ chain
and  $ (0.08610, 0.12669) $ for  $\xi_2$ chain.
With AMH-Full model, the $95\%$  CI for $\xi_1$ chain is $ (0.78294, 1.31152)$
and $\xi_2$ chain is $ (0.07024, 0.11360).$ }
\label{inv_hist101}
\end{figure}

\begin{figure}[bt]
\centering
\includegraphics[width=1.0\textwidth,height=0.550\textwidth]{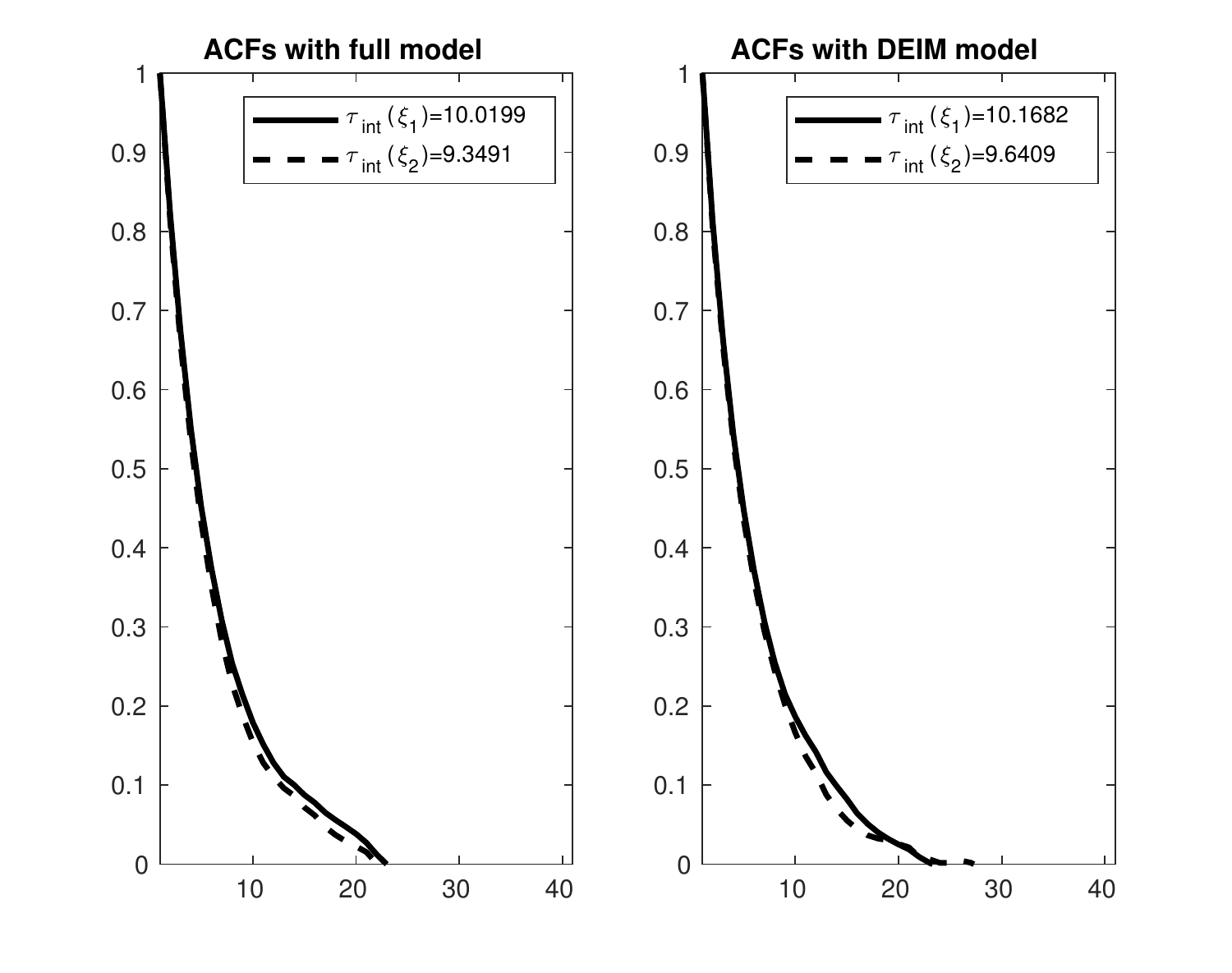}
\caption{Autocorrelation functions (ACFs)  for 
  $\xi_1$ and $\xi_2$ chains computed with AMH-Full (left) and  AMH-DEIM (right) models. These functions are  plotted against $J.$}
\label{inv_acf101}
\end{figure}

\begin{figure}[bt]
\centering
\includegraphics[width=0.9\textwidth,height=0.85\textwidth]{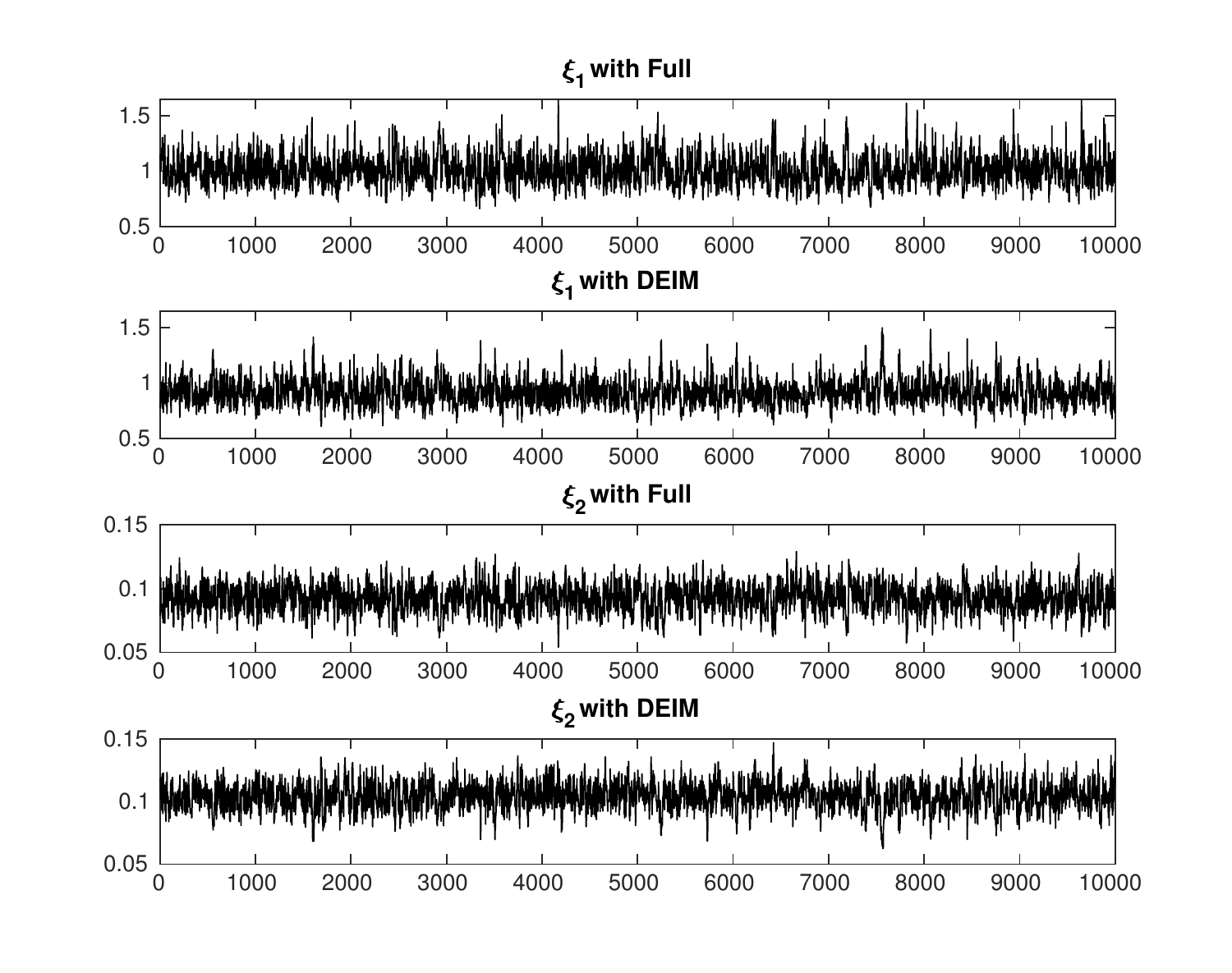}
\caption{MCMC samples
for the parameters $\xi  = [\xi_1, \xi_2]$ using AMH-Full (first and third) and  AMH-DEIM (second and fourth) models.
}
\label{inv_mc101}
\end{figure}

In addition to the above statistics,  the 
histograms, the scatter plots, as well as the pairwise plots of the chains from the
AMH-DEIM and AMH-full models are shown. 
Finally, we compare the computational times required by the AMH-DEIM and AMH-full models
to solve the statistical inverse problem.

\begin{figure}[bt]
\centering
\includegraphics[width=1.0\textwidth,height=0.50\textwidth]{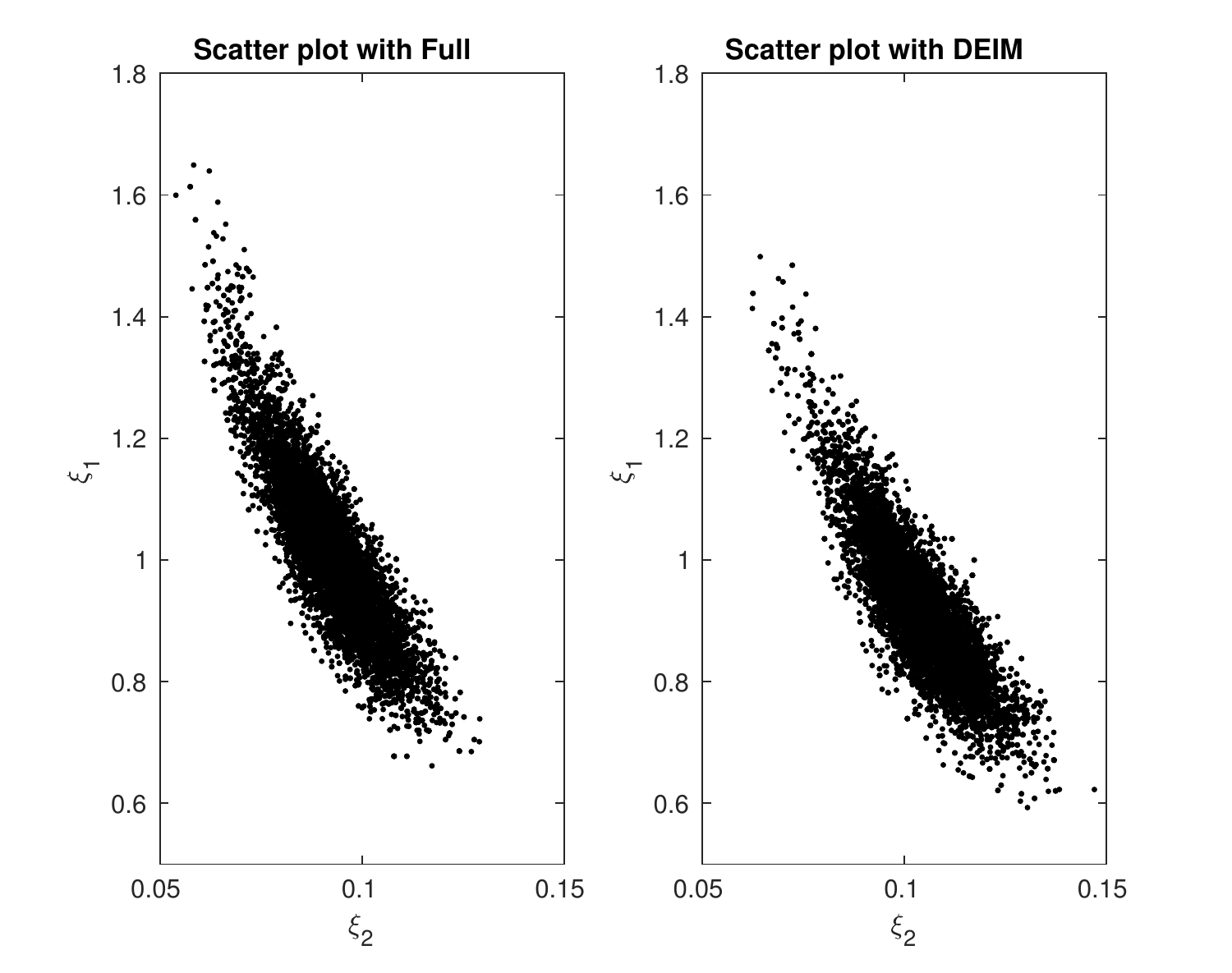}
\caption{Pairwise (scatter) plot of AMH-Full (left) and AMH-DEIM (right) samples of the parameters $\xi  = [\xi_1, \xi_2].$}
\label{inv_scatter101}
\end{figure}

The performance of the Markov chain from the AMH-DEIM model (with $k=n=100$) and AMH-Full (with $N=16384$) are reported
in  Figures
\ref{inv_hist101}, \ref{inv_acf101}, \ref{inv_mc101}, and \ref{inv_scatter101}.
First, we show the histograms for the  posterior distributions for the two parameters  $\xi_1$ and $\xi_2$   in Figure
\ref{inv_hist101}, under which we have also reported the 
the $95\%$  confidence intervals (CIs).  With AMH-DEIM (right), the CI is $ (0.75923, 1.18712)$ for $\xi_1$ chain
and  $ (0.08610, 0.12669) $ for  $\xi_2$ chain. The  midpoints of these two intervals are, respectively, 
$\xi_1= (0.75923 + 1.18712)/2= 0.973175$ and $ \xi_2= (0.08610+ 0.12669)/2=0.106395;$
thus, the mean
values $[0.973175, 0.106395]$ associated with the calculated distribution of $\xi$
can be seen as  AMH-DEIM model's estimate for the true parameters $\xi  = [1, 0.1].$
Indeed, the Geweke $p$-values obtained with  the AMH-DEIM algorithm for both $\xi_1$ and $\xi_2$ chains are  $0.99835$ 
and $0.98396,$ respectively. Observe that each of them   is over $98\%;$  this suggests that both chains are in equilibrium.

Figure \ref{inv_hist101} also shows the histograms for the  posterior distributions\footnote{Note 
that the histograms from AMH-Full
and AMH-DEIM are not exactly the same. The  difference is due
essentially to randomness and not the accuracy of the forward solutions.
If we ran the two models with the same seed (of random numbers), they would look identical.}
for both  parameters, 
with the AMH-Full model (left);
in this figure it is also reported that   
the $95\%$  CI for $\xi_1$ chain is $ (0.78294, 1.31152)$
and $\xi_2$ chain is $ (0.07024, 0.11360).$ Note that the midpoints of these intervals 
are, respectively,  
$1.0473$ and $ 0.09192.$ Thus, it can be seen that the estimates from the 
AMH-DEIM model (i.e., $0.973175, $ and $0.106395$) are better approximations of the true parameters
$\xi  = [1, 0.1].$
The Geweke $p$-values obtained with  the AMH-full model for both $\xi_1$ and $\xi_2$ chains are  $0.99808$ 
and $0.99441,$ respectively. Similarly to the AMH-DEIM case, these $p$-values indicate that 
the chains are in equilibrium.

Figure \ref{inv_acf101} shows the ACF for the two chains obtained with AMH-Full (left) and AMH-DEIM (right) schemes.
Note that  the  ACFs of $\xi_1$ and $\xi_2$ can be seen to approach zero smoothly. Observe also from the figure
that, with  the AMH-DEIM scheme, the IACT for the two chains  $\xi_1$ and $\xi_2$ are
 $10.16$ and $9.640,$  respectively. This suggests that 
 roughly every $10$th sample of both chains from AMH-DEIM model is independent.
 Similarly, with the AMH-Full scheme, the IACT for the    $\xi_1$ and $\xi_2$ chains are
 $10.019$ and $9.349,$  respectively.
 This suggests that 
 approximately every $10$th sample from $\xi_1$ chain and every $9$th sample from $\xi_2$ chains  are independent.
 Besides,
 the Markov chains plots in Figure \ref{inv_mc101}  and the scatter plots in Figure \ref{inv_scatter101},
together with the above statistical checks, imply that the computed Markov chains from both AMH-DEIM and AMH-Full schemes
``adequately'' approximate  the  posterior probability distributions for the unknown input parameters $\xi  = [\xi_1, \xi_2]$
in the inverse problem.

Undoubtedly, the computational complexity associated with the reduced models 
is significantly smaller than that of the full model. This 
justifies why we replace the full  nonlinear solver with the reduced solvers.
The overall consequence is then the reduction in the computational complexity of the solution of the statistical 
inverse problem   using the MCMC model (AMH). In general, 
in all our computations  the preconditioned
 AMH-DEIM model reduces the computational complexity of the statistical
inverse problem by $50\% - 80\%.$
More precisely, while  the AMH-Full model (with $N=16384$)
solves the inverse problem considered above in $5809.8\;  \mbox{sec}$  $ (96.3\;  \mbox{mins}),$  it
takes the AMH-DEIM  model with $k=n=100$ about $1400.16\;  \mbox{sec}$  $ (23.3\;  \mbox{mins})$,
which means the DEIM model can reduce the computational time by about $76\%.$ 

\bibliographystyle{siam}
\bibliography{akwumref}
\end{document}